\documentclass[a4paper, 11pt, onecolumn, oneside]{article}
\usepackage[T1]{fontenc}
\usepackage[utf8]{inputenc}
\usepackage{makeidx}
\usepackage{graphicx}
\usepackage{float}
\usepackage{amssymb}
\usepackage{amsmath}
\usepackage{latexsym}
\usepackage{wasysym}
\usepackage{amsthm}
\usepackage{natbib}
\usepackage{}
\usepackage{amsfonts}
\usepackage{amsmath}
\usepackage{amssymb}
\usepackage{amstext}
\newtheorem{theorem}{Theorem}
\newtheorem{proposition}[theorem]{Proposition}

\begin{document}

\title{Extreme Value  distribution for  singular measures}

\author{Faranda, Davide\\
\small{\textit{Department of Mathematics and Statistics, University of Reading;}}\\
\small{\textit{Whiteknights, PO Box 220, Reading RG6 6AX, UK.} d.faranda@pgr.reading.ac.uk}\\ \\
Lucarini, Valerio\\
\small{\textit{Department of Meteorology, University of Reading;}}\\
\small{\textit{Department of Mathematics and Statistics, University of Reading;}}\\
\small{\textit{Whiteknights, PO Box 220, Reading RG6 6AX, UK.} v.lucarini@reading.ac.uk}\\ \\
Turchetti, Giorgio\\
\small{\textit{Department of Physics, University of Bologna.INFN-Bologna}}\\
\small{\textit{Via Irnerio 46, Bologna, 40126, Italy.} turchett@bo.infn.it}\\ \\
Vaienti, Sandro\\
\small{\textit{UMR-6207, Centre de Physique Th\'eorique, CNRS, Universit\'es d'Aix-Marseille I,II,}}\\
\small{\textit{Universit\'e du Sud Toulon-Var and FRUMAM}}\\
\small{\textit{(F\'ed\'eration de Recherche des Unit\'es de Math\'ematiques de Marseille);}}\\
\small{\textit{CPT, Luminy, Case 907, 13288 Marseille Cedex 09, France.}}\\
\small{vaienti@cpt.univ-mrs.fr}\\
}
\date{}
%
%
%

\maketitle

\begin{abstract}

In this paper we perform an analytical and numerical study of Extreme Value distributions in discrete dynamical systems that have a singular measure. Using the block maxima approach described in \citet{faranda2011numerical} we show that, numerically, the Extreme Value distribution for these maps can be associated to the Generalised Extreme Value family  where the parameters scale with the information dimension. The numerical analysis are performed on a few low dimensional maps. For the middle third Cantor set and the Sierpinskij triangle obtained using Iterated Function Systems,  experimental parameters show a very good agreement with  the  theoretical values. For strange attractors like Lozi and H\`enon maps a slower convergence to the Generalised Extreme Value distribution is observed. Even in presence of large statistics the observed convergence is slower if compared with the maps which have an absolute continuous invariant measure. Nevertheless and within the uncertainty computed range, the results are in good agreement with the theoretical estimates.

\end{abstract}

\textbf{The existence of extreme value laws for dynamical  systems
preserving an absolutely continuous invariant measure or a singular
continuous invariant measure has been recently proven if strong
mixing properties or  exponential hitting time statistics on balls
are  satisfied.  In our previous work we have shown that there
exists an algorithmic way to study extrema by using a block-maxima
approach for dynamical systems which possess an absolutely
continuous invariant measure and satisfy certain mixing properties.
In this work we test our algorithm for maps that do not have an
absolutely continuous invariant measure showing that  the cumulative
distribution function of maxima is related to the scaling of the
measure of a ball centered around generic points. The scaling
exponent turns out to be the Hausdorff dimension of the measure
(also known as information dimension). Even if we cannot estimate
analytically the asymptotic behavior of  the measure of the balls,
the agreement with the numerical simulations we have carried out for
different maps suggests the validity of our proposed scaling in
terms of the information dimension. Our conjecture has been tested
with numerical experiments on different low dimensional maps such as
the middle third Cantor set, the Sierpinskij triangle, Iterated
Function System (IFS) with non-uniform weights and strange
attractors such as Lozi and H\'enon. In all cases considered, there
is a good agreement between the theoretical parameters and the
experimental ones although, in case of strange attractors which
exhibit multifractal structures, the convergence is slower. To
perform the numerical simulations it has been used the L-moments
procedure in order to overcome the difficulties of dealing with a
singular continuous invariant measure.}

\section{Introduction}
\subsection{Classical Extreme Value Theory}
Extreme Value Theory (EVT), developed  for the study of stochastical series of independent and identical distributed variables  by   \citet{fisher} and formalized by \citet{gnedenko}, has been successfully applied to different scientific fields to understand and possibly forecast events that occur with very small probability but that can be extremely relevant from an economic or social point of view:  extreme floods \citep{gumbel1941}, \citep{sv}, \citep{friederichs}, amounts of large insurance losses \citep{brodin}, \citep{cruz2002}; extreme earthquakes \citep{sornette}, \citep{cornell}, \citep{burton1979seismic}; meteorological  and climate events \citep{felici1},  \citep{vitolo},  \citep{alttresh}, \citep{nicholis1997clivar}, \citep{smith1989extreme}. An extensive review of the techniques and applications related to the EVT is presented in  \citet{ghil2010extreme}.\\

\citet{gnedenko}  studied the convergence of maxima of i.i.d.
variables $$X_0, X_1, .. X_{m-1}$$ with cumulative distribution
(cdf) $F(x)$ of the form:
$$F(x)=P\{a_m(M_m-b_m) \leq x\}$$
Where $a_m$ and $b_m$ are normalizing sequences and $M_m=\max\{ X_0,X_1, ..., X_{m-1}\}$. It may be rewritten as $F(u_m)=P\{M_m \leq u_m\}$ where $u_m=x/a_m +b_m$. Under general hypotesis on the nature of the parent distribution of data, \citet{gnedenko} show that the distribution of maxima, up to an affine change of variable,  obeys to one of the following three laws:
\begin{itemize}
\item Type 1 ({\em Gumbel}). \begin{equation}E(x)=\exp(-e^{-x}), \ -\infty<x<\infty \label{gum}\end{equation}
\item Type 2 ({\em Fr\'echet}). \begin{equation}E(x)= \begin{cases}
 0,\  x\le 0 \\
\exp(-x^{-\xi}),\  \mbox{for some} \ \xi>0, \ x>0
\end{cases} \label{fr} \end{equation}

\item Type 3 ({\em Weibull}). \begin{equation} E(x)= \begin{cases}
\exp(-(-x)^{-\xi}),\  \mbox{for some} \  \xi>0, \ x\le 0 \\
1,\  x>0
\end{cases} \label{wei}\end{equation}

\end{itemize}

Let us define the right endpoint $x_F$ of a distribution function $F(x)$ as:
\begin{equation}
x_F=\sup\{ x: F(x)<1\}
\end{equation}

then, it is possible to compute normalizing sequences   $a_m$ and $b_m$ using the following corollary of Gnedenko's theorem :\\
\textbf{Corollary (Gnedenko):}  \textit{The normalizing sequences $a_m$ and $b_m$ in the convergence of normalized maxima $P\{a_m(M_m - b_m) \leq x\} \to F(x)$ may be taken (in order of increasing complexity) as:}

\begin{itemize}

\item \textit{Type 1:} $\quad a_m=[G(\gamma_m)]^{-1}, \quad b_m=\gamma_m$;

\item \textit{Type 2:} $\quad a_m=\gamma_m^{-1}, \quad b_m=0$;

\item \textit{Type 3:} $\quad a_m=(x_F-\gamma_m)^{-1}, \quad b_m=x_F$;

\end{itemize}
\textit{where}
\begin{equation}
\gamma_m=F^{-1}(1-1/m)=\inf\{x; F(x) \geq 1-1/m\}
\label{gamma}
\end{equation}

\begin{equation}
G(t)=\int_t^{x_F} \frac{1-F(u)}{1-F(t)}du, \quad  t<x_F
\label{gneg}
\end{equation}

In \citep{faranda2011numerical} we have shown that this approach is equivalent  to fit unnormalized data directly  to a single family of generalized distribution called GEV distribution with cdf:

\begin{equation}
F_{G}(x; \mu, \sigma,
\xi')=\exp\left\{-\left[1+{\xi'}\left(\frac{x-\mu}{\sigma}\right)\right]^{-1/{\xi'}}\right\}
\label{cumul}
\end{equation}

which holds for $1+{\xi'}(x-\mu)/\sigma>0 $, using $\mu \in
\mathbb{R}$ (location parameter) and $\sigma>0$ (scale parameter) as
scaling constants in place of $b_m$, and $a_m$ \citep{pickands}, in
particular, in \citet{faranda2011numerical} we have shown that the
following relations hold:

$$\mu=b_m \qquad \sigma=\frac{1}{a_m}. $$

${\xi'} \in \mathbb{R}$ is the shape parameter also called the tail index: when ${\xi'} \to 0$, the distribution corresponds to a
Gumbel type ( Type 1 distribution).  When the index is positive, it corresponds to a Fr\'echet (Type 2 distribution); when the index is negative, it corresponds to a Weibull (Type 3 distribution).\\

To analyze the extreme value distribution in a series of data two main approaches can be applied: the Peak-over-threshold approach and the Block-Maxima approach. The former consists in looking at exceedance over high thresholds \citep{todorovic1970stochastic} and a Generalized Pareto distribution  is used for modeling data obtained as excesses over thresholds \citep{smith1984threshold}, \citep{davison1984modelling}, \citep{davison1990models}.\\
The so called block-maxima approach is widely used in climatological and financial applications since it represents a very natural way to look at extremes in fixed time intervals: it consists of dividing the data series of some observable into bins of equal length and selecting the maximum (or the minimum) value in each of them \citep{coles}, \citep{felici1}, \citep{katz}, \citep{katz2}, \citep{katz2005statistics}.\\

\subsection{Extreme Value Theory for dynamical systems}

As far as the classical EVT is concerned, we should restrict our domain of investigation to the output of stochastic processes. Obviously, it is of crucial relevance for both mathematical reason and for devising a framework to be used in applications, to understand under which circumstances the time series of observables of deterministic dynamical system can be treated using EVT.\\
 Empirical studies show that in some cases a dynamical observable obeys to the extreme value statistics even if the convergence is highly dependent on the kind of observable we choose \citep{vannitsem}, \citep{vitolo}, \citep{vitolo2009robust}. For example,  \citet{balakrishnan} and more recently \citet{nicolis} and \citet{haiman} have shown that for regular orbits of dynamical systems we don't expect to find convergence to EV distribution.\\
The first rigorous  mathematical approach to extreme value theory in dynamical systems goes back to the pioneer paper by P. Collet in 2001 \citep{collet2001statistics}. Important contributions have successively been given by  \citet{freitas2008}, \citet{freitas},\citet{freitas2010extremal} and by \citet{gupta2009extreme}. The starting point of all these investigations was
to associate to the stationary stochastic process given by the dynamical
 system, a new stationary independent sequence which enjoyed one of the
  classical three extreme value laws, and this laws could be pulled back
  to the original dynamical sequence. To be more precise we will consider
   a dynamical system $(\Omega, {\cal B}, \nu, f)$, where $\Omega$ is the
    invariant set in some manifold, usually $\mathbb{R}^d$, ${\cal B}$ is the
     Borel $\sigma$-algebra, $f:\Omega\rightarrow \Omega$ is a measurable map
      and $\nu$ a  probability $f$-invariant Borel measure. The stationary
       stochastic process given by the dynamical system will be of the form
       $X_m=g\circ f^m$, for any $m\in \mathbb{N}$, where the observable $g$
       has values in $\mathbb{R}\cup\pm\infty$ and achieves a global maximum
       at the point $\zeta\in \Omega$. We therefore study the partial maximum
       $M_m=\max\{X_0, \dots, X_{m-1}\}$, in particular we look for normalising
       real  sequences $\{a_m\}, \{b_m\}, m\in \mathbb{R}^{+}$ for
        which $\nu\{x; a_m(M_m-b_m)\le t\}=\nu\{x; M_m\le u_m\}$ converge
        to a non-degenerate distribution function; here $u_m=\frac{t}{a_m}+b_m$
        is such that $m\nu(X_0>u_m)\rightarrow \tau$, for some positive $\tau$
        depending eventually on $t$: we defer to the book \citep{leadbetter} for
         a clear and complete picture of this approach.  We will associate to
         our process a new i.i.d. sequence $\tilde{X}_0, \cdots, \tilde{X}_{m-1}$
         whose distribution is the same as that of $X_0$ and with partial maximum: $\tilde{M}_m=\max\{\tilde{X}_0, \cdots, \tilde{X}_{m-1}\}$. Properly normalized the distribution of such a maximum converges to one of  the three laws in equations \ref{gum}-\ref{wei} and this is the interesting content of the Extreme Value Theory. Equations \ref{gum}-\ref{wei}  will be satisfied by our original process too,
 whenever we would be able to prove that
$$
\lim_{m\rightarrow \infty}\nu(\tilde{M}_m\le u_m)=\lim_{m\rightarrow \infty}\nu(M_m\le u_m)
$$
This can be achieved if one can prove  two sufficient conditions called $D_2$ and $D'$ and which we briefly quote and explain in the footnote: these conditions  basically require a sort of independence of the stochastic dynamical sequence in terms of uniform mixing condition on the distribution functions. In particular condition $D_2$, introduced in its actual  form by Freitas-Freitas \citet{freitas2008},   could be checked directly by estimating the rate of decay of correlations for H\"older observables \footnote{We briefly state here the two conditions, we defer to the next section for more details about the quantities introduced. If $X_m, m\ge 0$ is our stochastic process, we can define  $M_{j,l}\equiv \{X_j, X_{j+1},\cdots,X_{j+l}\}$ and we put $M_{0,m}=M_m$.  The condition $D_2(u_m)$ holds for the sequence $X_m$ if for any integer $l, t,m$ we have $|\nu (X_0>u_m, M_{t,l}\le u_m)-\nu (X_0>u_m)\nu (M_{t,l}\le u_m)|\le \gamma(m,t)$, where $\gamma(m,t)$ is non-increasing in $t$ for each $m$ and $m\gamma(m,t_m)\rightarrow 0$ as $m\rightarrow \infty$ for some sequence $t_m=o(m)$, $t_m \rightarrow \infty$.\\ We say condition $D'( u_m)$ holds for the sequence $X_m$ if $\lim_{l\rightarrow \infty}\limsup_{m} m\sum_{j=1}^{[m/l]}\nu(X_0>u_m, X_j>u_m)=0$. }.  Another interesting issue of of the previous works was the choice of the observables $g$'s: it is chosen as  a function $g(\mbox{dist}(x, \zeta))$ of the distance with respect to a given point $\zeta$,   with the aim that $g$ achieves a global maximum at almost all points $\zeta\in \Omega$; for example $g(x)=-\log x$. In particular the observable $g$ was taken in one of three different classes  $g_1, g_2, g_3$, see Sect. 2 below, each one being again a function of the distance with respect to a given point $\zeta$. The choice of these particular forms for the $g$'s is just to fit with the necessary and sufficient condition on the tail of the distribution function $F(u)=\nu\{x; X_0\le u\}$,  in order to exist a non-degenerate limit distribution for the partial maxima \citep{freitas}, \citep{holland}. We use here the fact that, thanks to conditions $D_2$ and $D'$, the distributions of $X_0$ rules out the distribution of our non-independent process $X_m$ as well. It is important to remind that the previous conditions will determine the exponent $\xi$ in the types 2 and 3 for $E(x)$.  \\ Another major step in this field   was achieved by establishing a connection between the extreme value laws and the statistics of first return and hitting times, see the papers by \citet{freitas} and \citet{freitasNuovo}. They showed in particular that for dynamical systems preserving an absolutely continuous invariant measure or a singular continuous invariant measure $\nu$, the existence of an exponential hitting time statistics on balls around $\nu$ almost any point $\zeta$ implies the existence of extreme value laws for one of the observables of type $g_i, i=1,2,3$ described above. The converse is also true, namely if we have an extreme value law which applies to the observables of type $g_i, i=1,2,3$ achieving a maximum at $\zeta$, then we have exponential hitting time statistics to balls with center $\zeta$. Recently these results have been generalized to local returns around balls centered at periodic points \citep{freitas2010extremal}. \\
 In the context of singular measures, the EVT has been developed
in the recent paper by \citet{freitasNuovo}. The main goal of their
paper was to establish a connection with hitting and return time
statistics; for that purpose they considered returns in balls and
also into cylinders. We are particularly interested in their Theorem
1, about balls, since it covers the class of observables considered
in this paper; in particular we use here one direction of the
theorem which allows us to get the extreme value distributions if
the exponential  return time statistics has been previously
established for balls centered around almost all points and with
respect to the probability invariant measure (in this manner we do
not need to check conditions $D_2$ and $D'$). 

\subsection{This work}

In our previous work \citet{faranda2011numerical} we have shown that there exists an algorithmic way to study EVT by using a block-maxima approach for  dynamical systems which possess an absolutely continuous invariant measure and satisfy the mixing properties given by conditions $D_2$ and $D'$.  We have established the best conditions to observe convergence to the analytical results highlighting  deviations from theoretical expected behavior depending on the number of maxima and number of block-observation. Furthermore, we have verified that the normalising process of variables can be applied a posteriori and a fit of unnormalised data produce a distribution that belongs to the Generalised Extreme Value (GEV) distributions family.\\
In this work we test our algorithm for maps that do not have an absolutely continuous invariant measure. We remind that in the context of dynamical system, the invariant measure plays the role of the probability measure on the space of  events; in this respect the general theory of extremes will continue to apply no matter such a probability is absolutely continuous or singular with respect to Lebesgue. The interesting point is that for the choice of observables we did (the functions $g_i$), the cumulative distribution function $F$  will be related to the scaling of the measure $\nu(B_r(z))$ of a ball  $B_r(z)$ of radius $r$ and centered at the  point $z$, and such a scaling exponent turns out to be the Hausdorff dimension of the measure (also known as information dimension), when the point $z$ is generically chosen. The experimental and accessible  parameters of the GEV distributions will be explicitly expressed in terms of such a dimension.\\
In order to get the values of
  $\xi$ and of $a_m$ and $b_m$ for finite $m$
  one should know how the measure   of the ball $B_r(z)$
   behaves as a function of $r$ and of $z$
    and for measures which are not absolutely continuous.
    We notice that for absolutely continuous measure that
    approach works, at least in a few cases, and we quote
    our previous paper for that. Instead for singular continuous measures
    like those supported on Cantor sets,
    we are not aware of any analytic result allowing
     to get the few orders expansion of $\nu(B_r(z))$.
      This will prevent us to compute rigorously the
      normalising constant $a_m$ for type 1 observables $g_1$;
       instead  we will get the the limiting values of $b_m$
       for type 1 and the limiting values of $a_m$  for type 2 and 3.
       Moreover we could not compute rigorously the exponent $\xi$.
       The values proposed for those non-rigorous constants are obtained
        by simply  approximating $\nu(B_r(z))$ with $r^{D}$.
         The agreement with the numerical simulations suggests
         that there were good choices and suggests also a direct
         proof of the EVT for our observables and with the normalising
         constants indicated by our heuristic analysis.\\
As explained in the previous subsection we can either check the conditions $D_2$ and $D'$ or the existance of an exponential return time statistics. The latter is the case of iterated function systems considered in section 3.2: these are in fact given by expanding maps (since they verify
the so-called {\em open set condition}) and the exponential return
time statistics for balls could be proved, for instance, using the
technique in \citet{bessis1987mellin}. It will be also the case for
the H\'enon attractor with the parameters studied by Benedicks and
Carlesson: for those parameters the attractor exists and carries an
SRB measure; moreover very recently Chazottes and Collet established
the Poissonian statistics for the number of visits in balls around
generic points w.r.t. the SRB measure. Our numerical computation
will concern instead the {\em usual} H\'enon attractor. Finally, we
will consider the Lozi attractor, and in this case we will quote the
result by \citet{gupta2009extreme}, which proves the existence of
the extreme value distributions for the observables constructed with
the functions $g_i$ and for balls around almost any point w.r.t. the
SRB measure.  As a final remark, we stress that the results by
Freitas-Freitas and Todd have been proved under the assumption that
$\nu(B_r(\zeta))$ is a continuous function of $r$: this is surely
true for all the previous examples and such a condition will play a
major role in our next considerations too.\\
This work is organized as follows: in Section 2 we present the analytical results for the EVT in maps with singular measures deriving the asymptotic behavior of normalising sequences and parameters. In Section 3 we present the numerical procedure used for the statistical inference of the GEV distribution and the numerical experiments that we have carried out for both singular measures generated with Iterated Function Systems and maps with a less trivial measures such as the Baker transformation, H\`enon and Lozi maps. Eventually, in Section 4 we present our conclusion and proposal for future work.

\section{Extreme Value Theory for maps with singular measures}
\label{background}

\subsection{Definitions and Remarks}
Let us consider a dynamical systems $(\Omega, {\cal B}, \nu, f)$, where $\Omega$ is the
    invariant set in some manifold, usually $\mathbb{R}^d$, ${\cal B}$ is the
     Borel $\sigma$-algebra, $f:\Omega\rightarrow \Omega$ is a measurable map
      and $\nu$ a  probability $f$-invariant Borel measure.\\
As we said in the Introduction and in order to adapt the extreme value theory to dynamical systems, we will consider the stationary stochastic process $X_0,X_1,...$  given by:

\begin{equation}
X_m(x)=g(\mbox{dist}(f^m (x), \zeta)) \qquad \forall m \in \mathbb{N}
\label{sss}
\end{equation}

where 'dist' is a distance on the ambient space  $\Omega$, $\zeta$ is a given point and $g$ is an observable function, and whose partial maximum is defined as:

\begin{equation}
{M_m}= \max\{ X_0, ... , X_{m-1} \}
\label{maxi}
\end{equation}

The probability measure will be here the invariant measure $\nu$ for the dynamical system. We will also suppose that our systems which verify the condition $D_2$ and $D'$ which will allow us to use the EVT for i.i.d. sequences. As we said above, we will use   three types of observables $g_i,i=1,2,3$,  suitable to obtain one of the three types of EV distribution  for normalized maxima:

\begin{equation}
g_1(x)= -\log(\mbox{dist}(x,\zeta))
\label{g1}
\end{equation}

\begin{equation}
g_2(x)=\mbox{dist}(x, \zeta)^{-1/\alpha}
\label{g2}
\end{equation}

\begin{equation}
g_3(x)=C - \mbox{dist}(x,\zeta)^{1/\alpha}
\label{g3}
\end{equation}

where $C$ is a constant and $\alpha>0 \in \mathbb{R}$.\\
These three type of functions are representative of broader classes which are defined, for instance, throughout equations (1.11) to (1.13) in Freitas et al. [2009]; we now explain the reasons and the meaning of these choices. First of all these functions have in common the following properties: (i) they are defined on the positive semi-axis $[0,\infty]$ with values into $\mathbb{R}\cup \{+\infty\}$; (ii) $0$ is a global maximum, possibly equal to $+\infty$; (iii) they are a strictly decreasing bijection in a neighborhood $V$ of $0$ with image $W$. Then we consider three types of behavior which generalize the previous specific choices:\\
{\em Type 1}: there is a strictly positive function $p:W\rightarrow  \mathbb{R}$ such that $\forall y\in \mathbb{R}$ we have $$ \lim_{s\rightarrow g_1(0)}\frac{g_1^{-1}(s+yp(s))}{g_1^{-1}(s)}=e^{-y}
$$
{\em Type 2}: $g_2(0)=+\infty$ and there exists $\beta>0$ such that $\forall y>0$ we have $$ \lim_{s\rightarrow \infty}\frac{g_2^{-1}(sy)}{g_2^{-1}(s)}=y^{-\beta}
$$
{\em Type 3}: $g_3(0)=D<+\infty$ and there exists $\gamma>0$ such that $\forall y>0$ we have $$ \lim_{s\rightarrow 0}\frac{g_3^{-1}(D-sy)}{g_3^{-1}(D-s)}=y^{\gamma}
$$

The Gnedenko corollary says that the different kinds of extreme value laws are determined by the distribution of 
\begin{equation}
F(u)=\nu (X_0\le u)
\label{fu}
\end{equation}

and by the right endpoint of $F$, $x_F$.\\

We need to compute and to control the measure $\nu( B_{r}(\zeta))$ of a ball of radius $r$ around the point $\zeta$. At this regard we will invoke, and assume, the existence of the following limit
\begin{equation}\label{MB}
\lim_{r\rightarrow 0}\frac{\log \nu( B_{r}(\zeta))}{\log r}, \ \mbox{for} \  \zeta \ \mbox{chosen} \ \nu-\mbox{a.e.}
\end{equation}
Moreover we will assume that $\nu( B_{r}(\zeta))$ is a continuous function of $r$ (see \citet{freitas} for a discussion of this condition which shows that all the examples considered in our paper will fit it).
When the limit (\ref{MB}) exists on a metric space equipped with the Borel $\sigma$-algebra and a probability measure $\nu$, it gives the {\em Hausdorff dimension of the measure} or {\em information dimension}, defined as the infimum of the Hausdorff dimension taken over all  the set of $\nu$ measure $1$ \citep{young1982dimension}. This limit could be proved to exist for a large class of dynamical systems and whenever $\nu$ is an invariant measure: let us indicate it with $\Delta$ without written explicitly its dependence on $\nu$. For example, for a very general class of
 one-dimensional maps with positive metric entropy, $\Delta$ is equal to the ratio between the metric entropy and the (positive) Lyapunov exponent of $\nu$ \citep{ledrappier1981some}. For two dimensional hyperbolic diffeomorphisms, $\Delta$ is equal to the product of the metric entropy times the difference of the reciprocal of the positive and of the negative Lyapunov exponents \citep{young1982dimension}. The information dimension is a lower bound of the Hausdorff dimension of the support of the measure $\nu$ and it is an upper bound of the correlation dimension \citep{yakov1998pesin}, \citep{hentschel1983infinite}, \citep{grassberger1983generalized}, \citep{bessis1988generalized}, \citep{bessis1987mellin}, \citep{cutler1989estimation}.\\

\subsection{Limiting behavior of the Extreme Value Theory parameters}

We summarize the three basic assumptions for the next considerations:
\begin{itemize}
\item Assumption 1: our dynamical system verifies conditions $(D_2)$ and $D'$.
\item Assumption 2: the measure of a ball is a continuous function of the radius for almost all the center points; moreover such a measure has no atoms.
\item Assumption 3: the limit (\ref{MB}) exists (and its value is called $\Delta$) at almost all points $\zeta$.
\end{itemize}

Equipped with these conditions it is now possible to compute rigorously a
few of the  expected parameters for the three types of observables.

\paragraph{ Case 1: $\mathbf{g_1}$(x)= -log(dist(x,$\mathbf{\zeta}$)).}

Substituting equation \ref{sss} into equation \ref{fu} we obtain that:
\begin{equation}
\begin{split}
1-F(u) &=1-\nu(g(\mbox{dist}(x, \zeta))\leq u) \\
      &= 1 - \nu( -\log( \mbox{dist}(x, \zeta)) \leq u) \\
      &= \nu( \mbox{dist}(x, \zeta) < e^{-u})=\nu(B_{e^{-u}}(\zeta))\\
\end{split}
\label{g1_e1}
\end{equation}
$$x_F = \sup\{ u ; F(u) < 1\}$$
To use Gnedenko corollary it is necessary to calculate $x_F$;
in this case $x_F=+ \infty$ as we will explain in the proof below.

According to Corollary 1.6.3 in \citet{leadbetter} for type 1
$a_m=[G(\gamma_m)]^{-1}$ and $b_m=\gamma_m=F^{-1}(1-\frac{1}{m})$.
 We now show how to get the limiting value of $\gamma_m$; a similar proof will hold for type II and III.

\begin{proposition}\label{PP}
Let us suppose that our system verifies Assumptions 1,2,3 above and let us consider the observable $g_1$; then:
$$
\lim_{m\rightarrow \infty} \frac{\log m}{\gamma_m}=\Delta
$$
\end{proposition}
{\em Proof}\\
 By our choice of the observable we have:  $1-F(\gamma_m)=\nu(B_{e^{-\gamma_m}}(\zeta))=\frac{1}{m}$; since the measure is not atomic and it varies continuously with the radius, we have necessarily that $\gamma_m\rightarrow \infty$ when $m\rightarrow \infty$.  Now we fix $\delta>0$ and small enough; there will be $m_{\delta, \zeta}$ depending on $\delta$ and on $\zeta$, such that for any $m\ge m_{\delta, \zeta}$ we have
\begin{equation}\label{B2}
-\delta\gamma_m\le \log \nu(B_{e^{-\gamma_m}}(\zeta))+\Delta\gamma_m\le \delta \gamma_m
\end{equation}
Since $\log m-\Delta\gamma_m=-[\log \nu(B_{e^{-\gamma_m}}(\zeta))+\Delta\gamma_m]$ and by using the bounds (\ref{B2}) we immediately have
$$
-\delta\gamma_m \le \log m-\Delta\gamma_m\le \delta\gamma_m
$$
which proves the Proposition. \\

It should be clear that the previous proposition will not give us the value of $\gamma_m$ and of $b_m$, which is equal to $\gamma_m$ for type I observables. We have instead a rigourous limiting behavior and we will pose in the following:
$$
\gamma_m=b_m\sim \frac{1}{\Delta}\log m
$$
The values for finite $m$ could be obtained if one would dispose of
the functional dependence of $\nu(B_r(\zeta)$ on the radius $r$ and
the center $\zeta$: this  has been achieved for non-trivial
absolutely continuous invariant measure in our previous paper \citep{faranda2011numerical}. The
same reason prevent us to get a rigorous limiting behavior for $a_m=[G(\gamma_m)]^{-1}$.
 The only rigorous statement we can do is
that $G(\gamma_m)=o(\gamma_m)$; this follows by adapting the
previous proof of the proposition to another result (see \citet{leadbetter})
which says that for type I observables one has $\lim_{m\rightarrow
\infty}n(1-F\{\gamma_m+xG(\gamma_m)\})=e^{-x}$, for all real $x$:
choosing $x=1$ gives us the previous domination result. In the
following and again for numerical purposes we will take
$$
a_m=[G(\gamma_m)]^{-1}\sim \frac{1}{\Delta}
$$
This follows easily by replacing in formula (\ref{gneg}) $\nu(B_r(\zeta))\sim r^\Delta$ for $r$ small. We finish this part by stressing that for our observable we expect $\xi=\xi'=0$.
\paragraph{Case 2: $\mathbf{g_2}$(x)=dist(x,$\mathbf{\zeta)^{-1/\alpha}}$.}

In this case we have

\begin{equation}
\begin{split}
1-F(u)&= 1 - \nu(\mbox{dist}(x, \zeta)^{-1/\alpha} \leq u) \\
      &= 1 - \nu( \mbox{dist}(x, \zeta) \geq u^{-\alpha})\\
      &=  \nu( B_{u^{-\alpha}}(\zeta))
\end{split}
\end{equation}
and $x_F=+ \infty$. Since $b_m=0$ we have only to compute $a_m$ which is the reciprocal of $\gamma_m$ which is in turn defined by $\gamma_m=F^{-1}(1-1/m)$. By adapting Proposition \ref{PP} we immediately get that
$$
\lim_{m\rightarrow \infty}\frac{\log m}{\log \gamma_m}=\alpha \Delta
$$
which we allow us to use the approximation $a_m\sim
\frac{1}{m^{\frac{1}{\alpha \Delta}}}$ The exponent $\xi$ for Type II
observables is given by the following limit (see \citep{leadbetter},
Th. 1.6.2)
$$
\lim_{t\rightarrow \infty}(1-F(tx))/(1-F(t))=x^{-\xi}, \ \xi>0, \
x>0
$$
The crude approximation $\nu(B_r(\zeta))\sim r^\Delta$ for $r$ small,
will give immediately that $\xi\sim \alpha \Delta$ and this value will
appear in the exponent $\xi'=1/\xi$ in the distribution function
given by the GEV.

\paragraph{Case 3: $\mathbf{g_3}$(x)=C-dist(x,$\mathbf{\zeta)^{1/\alpha}}$.}

We have first of all:

\begin{equation}
\begin{split}
1-F(u)&= 1 - \nu(C- \mbox{dist}(x, \zeta)^{1/\alpha} \leq u) \\
      &=  \nu( B_{(C - u)^{\alpha}}(\zeta)) \\
\end{split}
\end{equation}

In this case $x_F= C<\infty$ and $a_m=(C-\gamma_m)^{-1}$; $b_m=C$.
The previous proposition immediately  shows that $\lim_{n\rightarrow
\infty}\frac{\log m}{-\alpha \log (C-\gamma_m)}=\Delta$ which gives the
asymptotic scaling $\gamma_m\sim C-\frac{1}{m^{\frac{1}{\alpha
\Delta}}}$; $a_m\sim m^{\frac{1}{\alpha \Delta}}$; $b_m=C$. Finally the
exponent $\xi$ is given again by Th. 1.6.2 in \citet{leadbetter} by
the formula
$$
\lim_{h\rightarrow 0}(1-F(C-hx))/(1-F(C-h))=x^{\xi}, \ \xi>0, \ x>0
$$
which with our usual approximation furnishes $\xi \sim \alpha \Delta$.

\section{Numerical Experiments}

\subsection{Procedure for statistical inference}
For the numerical experiments we have used a wide class of maps that
have singular measure, also considering the case of strange
attractors such as the ones observed  by iterating Lozi
 Map or H\`enon map. The algorithm used
is the same described in \citet{faranda2011numerical}: for each map
we run a long simulation up to $k$ iterations starting from a given
initial condition. From the trajectory we compute the sequence of
observables $g_1, g_2, g_3$ dividing it into $n$
bins each containing $m = k/n$ observations and eventually obtaining the empirical cdf of maxima.\\
In \citet{faranda2011numerical}  we have
used a Maximum Likelihood Estimation   (MLE) procedure working both on pdf and cdf (cumulative distribution
function), since our distributions were absolutely continuous and
the minimization procedure was well defined. In this case, we don't
have anymore the pdf  and consequently the fitting procedure via MLE
could give us wrong results. To avoid these problems
  we have used an L-moments estimation
   as detailed in \citet{hosking1990moments}.
   This procedure is completely discrete and can be used both for
   absolutely  continuous or  singular continuous  cdf. The L-moments
   are summary statistics for probability distributions and data samples.
    They are analogous to ordinary moments which meant that
    they provide measures of location, dispersion, skewness,
     kurtosis, but are computed from linear combinations of
     the data values, arranged in increasing order (hence the prefix L). Asymptotic approximations to sampling distributions are better for L-moments than for ordinary moments [\citet{hosking1990moments}, Figure 4].  The relationship between the moments and the parameters of the GEV distribution are described in \citet{hosking1990moments}, while the 95\% confidence intervals has been derived using a bootstrap procedure. As comparison we have checked that the results presented in \citet{faranda2011numerical} are comparable with L-moments methods.   We have found that both the methods give similar results even if L-moments has an uncertainty on the estimation of parameters  generally slightly bigger.\\
The empirical cdf contains plateaux which correspond to non accessible distances in correspondence of the holes of the Cantor set. The discrete nature of L-moments allows to overcome difficulties that may arise in singular continuous cdf: the normalization procedure carried out with this method consists in dividing each quantity computed via  L-moments by a function of the total number of data will prevent us from obtain a unnormalised distribution.  The last issue we want to address is the choice of a suitable model for our data: in principle, using L-moments procedure we can fit the data to any kind of known cdf. To validate the use of the GEV model we proceeded in the following way:\\
-\textit{A priori} the choice of a  GEV model arises naturally if the assumptions presented in the section 3 are satisfied. In this set up we can directly compute the parameters of GEV distribution using L-moments as described in \citet{hosking1990moments}.\\
-\textit{A posteriori} we can verify the goodness of fit to GEV family if we apply some  parametric or non-parametric tests commonly used in statistical inference procedures.  For this purpose we have fitted our experimental data to a wide class of well known  continuous distributions. Using Kolmogorov Smirnov test (see \citet{lilliefors} for a description of the test) we have measured the deviation between the empirical cdf and the fitted cdf, finding that  using the GEV distribution we effectively achieve a minimization of the deviation parameter. \\
We summarize below the results we expect from numerical experiments in
 respect to $n$ according to the conjecture described in the previous section.
 Since we keep the length of the series $k=n\cdot m$ fixed, the following
 relationships can be obtained simply replacing $m=k/n$ in the equations derived in the previous section.\\
For $g_1$ type observable:

\begin{equation}
\sigma= \frac{1}{\Delta} \qquad \mu \sim \frac{1}{\Delta}\ln(k/n) \qquad \xi'=0
\label{g1res}
\end{equation}

For $g_2$ type observable, we can either choose $b_m=0$ or $b_m=c\cdot m^{-\xi'}$ where $c \in \mathbb{R}$ is positive constant, as detailed in \citet{beirlant}. A priori, we do not know which
asymptotic sequences will correspond to the parameters $\mu$ in the experimental set up. The experimental procedure we use automatically select  $b_m=c\cdot m^{-\xi'}$, therefore the following results are presented taking into account  this asymptotic sequence:

\begin{equation}
\sigma\sim n^{-1/(\alpha \Delta)} \qquad \mu \sim n^{-1/(\alpha \Delta)} \qquad \xi'=\frac{1}{\alpha \Delta}
\label{g2res}
\end{equation}

For $g_3$ type observable:

\begin{equation}
\sigma\sim n^{1/(\alpha \Delta)} \qquad \mu = C \qquad \xi'=\frac{1}{\alpha \Delta}
\label{g3res}
\end{equation}

\subsection{IFS for Cantor Sets }
\label{CS}

A Cantor set can be obtained as an attractor of some Iterated Function Systems (IFS).  An IFS is a finite family of contractive maps $\{ f_1, f_2, ..., f_s \}$ acting on a compact normed space $\Omega$ with norm $|\cdot|$ and  possessing a unique compact limit  set (the attractor) $K \in \Omega$ which is non-empty
and invariant by the IFS, namely:

$$K = \bigcup\limits_{i=1}^{s} f_i(K).$$
We will put a few restrictions on the IFS in order to see it as the
 {\em inverse} of a genuine dynamical system; we will explain in a moment
 why this change of perspective will help us to compute observables on fractal
 sets. We defer to the fundamental paper by \citet{barnsley1985iterated},
  for all the material we are going to use.  \\ First of all we will consider
  the $f_i$ as strict contractions, namely there will be a number $0<\lambda<1$
  such that for all $i=1,\cdots,s$ we have $|f_i(x)-f_i(y)|<\lambda |x-y|$,
   for all $x,y\in X$.\\ Then we will suppose that each $f_i$ is one-to-one
   on the attractor $K$ and moreover $\forall i=1,\cdots,s$  we have
   $f_i(K)\cap f_j(K)=\emptyset$, $i\neq j$ (\textit{open set condition}). This will allow us to define
    a measurable map $T:K\rightarrow K$ by $T(x)=f_i^{-1}(x)$
     for $x\in f_i(K)$: the attractor $K$ will be the invariant
      set for the transformation $T$ which will play therefore
      the role of a usual dynamical system. A complete statistical
       description of a dynamical system is given by endowing it with
        an invariant probability measure; in particular we ask that
         this measure by ergodic if we want to compute the maxima of
          the sequence of events constructed with the observables
          $g_i$; we remind that these events are nothing but the
          evaluations along the forward orbit of an initial point chosen according to
          the measure. If, as always happens, the attractor $K$ has a
          fractal structure and zero Lebesgue measure, we could not
          get such an initial point for numerical purposes. We overcome this
           situation for attractors of global diffeomorphisms by
           taking the initial point in the basin of attraction and
           by
           iterating it: the orbit will be distributed according to
           the SRB measure (we will return later on this measure).
           For our actual attractors generated by non-invertible
           maps, the iteration of any point in the complement of the
           attractor will push the point far from it: it would be better to call
           repellers our invariant sets instead of attractors.  We have
           therefore to proceed in a different manner.
           The measures supported on the attractor $K$ will give the
           solution. First of all let us associated to each map
           $f_i$ a positive weight $p_i$ in such a way that
           $\sum_{i=1}^{s} p_i=1$. Then it is possible to prove the
           existence of a unique measure $\nu$ (called {\em
           balanced}) which enjoy the following properties:
           \begin{itemize}
           \item The measure $\nu$ is supported on the attractor $K$
           and it will be invariant for the map $T$ associated to
           our IFS (see above).
           \item For any measurable set $B$ in $X$ we have
           $$
           \nu(B)=\sum_{i=1}^{s}p_i\nu(f_i^{-1}(B))
           $$
           \item Let us put $(Sg)(x)= \sum_{i=1}^{s} p_i g(f_i(x))$, for
            a continuous functions $g$ on $\Omega$ and for
           {\em any} point $x\in \Omega$, then we have
           $$
           \lim_{n\rightarrow \infty} S^ng(x)= \int_\Omega g d\nu
           $$
           \end{itemize}
           This last item is very important for us; first of all it
           holds also for the characteristic function of a  set provided
           the boundary of this set has $\nu$ measure zero; therefore it
           is a sort of ergodic theorem because it states that the
           backward orbit constructed by applying to any point in $X$ the
           maps $f_i$ with weights $p_i$ will distribute on the
           attractor $K$ as the forward orbit (namely the orbit generated by the
           transformation $T$ associated to the IFS)
            of a point $y\in K$ and chosen
           almost everywhere according to $\nu$. \\ Let us give our
           first example.
           \subsubsection{Uniform weights and the Sierpisnkij
           triangle}
We consider  the middle one third Cantor set that is the attractor
of the IFS $\{ f_1, f_2 \}$ defined as:

\begin{equation}
\begin{cases}
f_1(x)= x/3 \mbox{ with weight } p_1\\
f_2(x)=(x +2)/3 \mbox{ with weight } p_2\\
\end{cases}
\label{IFS_C}
\end{equation}

where $x \in [0,1]$ and we set $p_1=p_2=1/2$
so that, at each time step, we have the same probability to iterate $f_1(x)$ or $f_2(x)$: 

 Equivalently the previous IFS can be written as:

\begin{equation}
x_{t+1}= (x_t + b)/3 \label{IFS_C2}
\end{equation}

where, at each time step,
we extract randomly with equal probability $b$ to be  0 or 2.\\

We will consider also the so called Sierpisnkij triangle,
defined by

\begin{equation}
\begin{cases}
x_{t+1}= (x_t + v_{p,1})/2 \\
y_{t+1}= (y_t + v_{p,2})/2 \\
\end{cases}
\label{IFS_S}
\end{equation}

We extract randomly at each time step and with equal probability the number
$p$ to be the integer 1,2 or 3. Then we iterate the map \ref{IFS_S} substituting the elements $v_{p,1}$ and $v_{p,2}$  of the following matrix:

\[ v = \left| \begin{array}{ccc}
1  & 0 \\
-1 & 0 \\
0 & 1 \end{array} \right|.\]

\begin{figure}[H] \begin{center}
\advance\leftskip-1.5cm
        \includegraphics[width=0.8\textwidth]{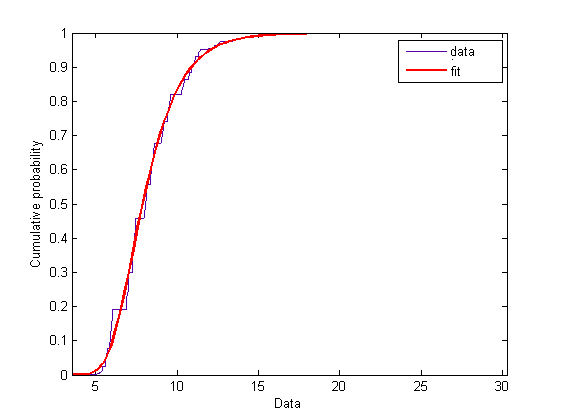}
        \caption{Empirical (blue) and fitted (red) cdf for IFS that generates a Cantor Set, $\zeta=0.775$, $g_1$ observable function.}
      \label{his_as}
    \end{center} \end{figure}

For these attracting sets the information dimensions are
well known, they are  $\Delta=\log(2)/\log(3)$ for the Cantor set and $\Delta=\log(3)/\log(2)$ for the Sierpinskij triangle \citep{sprott2003chaos}.\\
In the following experiments In order to choose the centers $\zeta$ of our balls we proceed by using again the backward preimages of any point in $\Omega$, namely we take a point $x\in \Omega$ and we consider $\zeta$ as  one of the preimages $f^{-t}(x)$ with $t$ much larger than the sequence of observed events; by what we said above, that preimage will be closer and closer to the invariant Cantor set and also it will approaches a generic point with respect to the balanced measure $\nu$.\\
First of all we have analised the  empirical cdf $F(u)$
of the extrema for $g_1$ observable. An example is shown  for the IFS in eq. \ref{IFS_C} in Figure \ref{his_as}.
The histogram is obtained iterating the map in equation \ref{IFS_C} for
  $5\cdot10^7$ iterations, $\zeta\simeq 0.775$, $\alpha=4$, $C=10$.
  Once computed the functions $g_1$ the series of maxima for each observable
  is computed taking each of them in bins containing 5000 values of $g_1$ for a total of 1000 maxima.
As claimed in the previous section, the cdf is a singular continuous
function and this is due to the structure of the Cantor set. The results are similar for  the other observables and other initial conditions.\\

To check that effectively the parameters of GEV distribution obtained by L-moments estimation are related to the fractal dimension of the attracting Cantor Set and the Sierpinskij triangle, we have considered an ensemble of $10^4$ different realizations of the eq. \ref{IFS_C} and eq. \ref{IFS_S}, starting from the same initial conditions. To check the behavior we have varied $n$ and $m$ keeping fixed the length of the series $k=10^7$. In \citet{faranda2011numerical} we have shown that a good convergence is observed when $n,m>1000$, therefore we will make all the considerations for $(n,m)$ pairs that satisfy this condition.\\

\begin{figure}[H] \begin{center}
\advance\leftskip-1.5cm
        \includegraphics[width=0.9\textwidth]{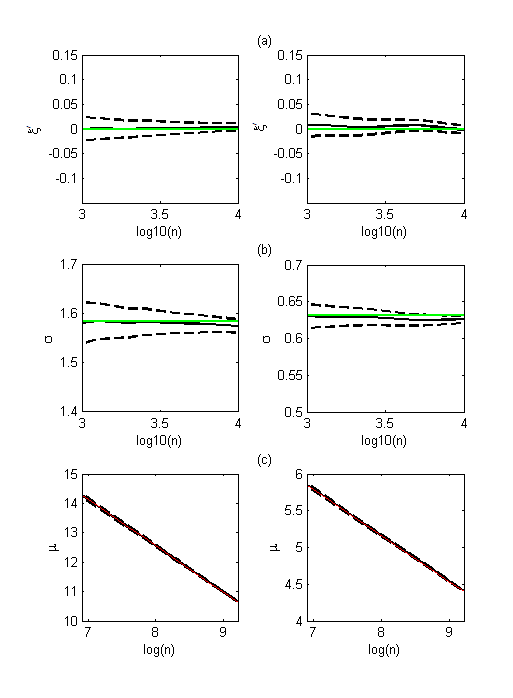}
        \caption{ $g_1$ observable. \textbf{a)} $\xi'$ VS $\log_{10}(n)$; \textbf{b)} $\sigma$ VS $\log_{10}(n)$; \textbf{c)} $\mu$ VS $\log(n)$.  Cantor set, Right: Sierpinskij triangle. Dotted lines represent one standard deviation,  red lines represent a linear fit, green lines are theoretical values.}
      \label{c1s1}
    \end{center} \end{figure}

\begin{figure}[H] \begin{center}
\advance\leftskip-1.5cm
        \includegraphics[width=0.9\textwidth]{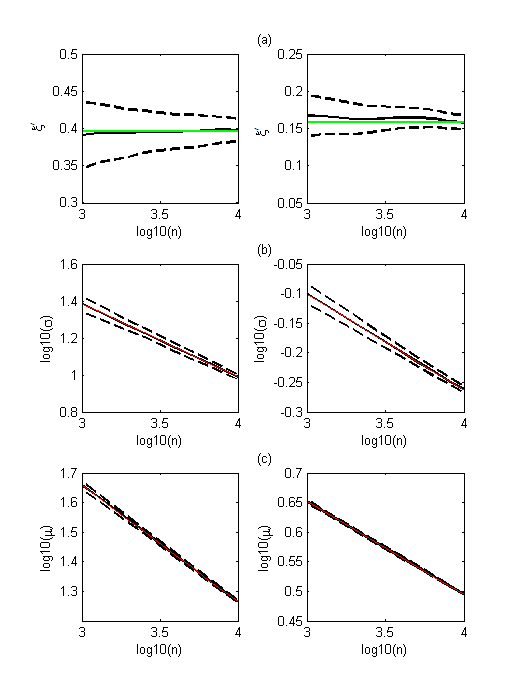}
        \caption{  $g_2$ observable \textbf{a)} $\xi'$ VS $\log_{10}(n)$; \textbf{b)} $\log_{10}(\sigma)$ VS $\log_{10}(n)$; \textbf{c)} $\log_{10}(\mu)$ VS $\log_{10}(n)$. Left: Cantor set, Right: Sierpinskij triangle. Dotted lines represent one standard deviation,
 red lines represent a linear fit, green lines are theoretical values.}
      \label{c2s2}
    \end{center} \end{figure}

\begin{figure}[H] \begin{center}
\advance\leftskip-1.5cm
        \includegraphics[width=0.9\textwidth]{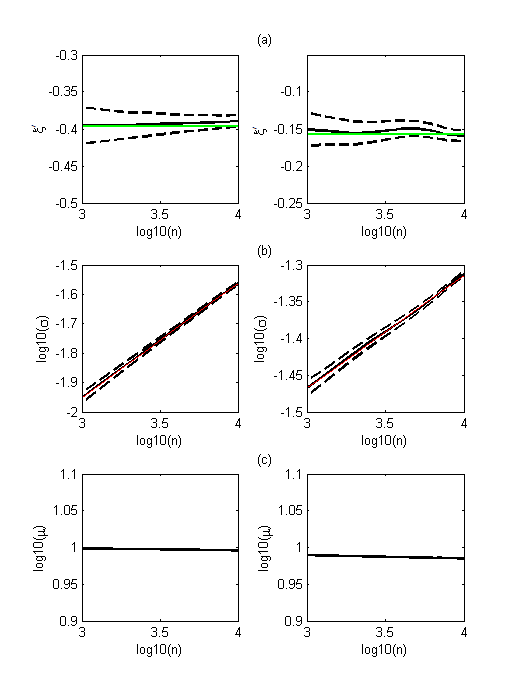}
        \caption{  $g_3$ observable. \textbf{a)} $\xi'$ VS $\log_{10}(n)$; \textbf{b)} $\log_{10}(\sigma)$ VS $\log_{10}(n)$; \textbf{c)} $\log_{10}(\mu)$ VS $\log_{10}(n)$.  Cantor set, Right: Sierpinskij triangle. Dotted lines represent one standard deviation,  red lines represent a linear fit, green lines are theoretical values.}
      \label{c3s3}
    \end{center} \end{figure}

In figures  \ref{c1s1}-\ref{c3s3} the results of the computation for the IFS that generates the Cantor Set ( plots on the left) and the Sierpisnkij triangle (plots on the right) are presented. In all the cases considered the behavior well reproduce the theoretical expected trend described in equations \ref{g1res}-\ref{g3res}. The initial condition here shown is $\zeta=1/3$ for Cantor Set and $\zeta\simeq(0.02, 0.40)$ for Sierpinskij triangle, but similar results hold for different initial conditions if chosen on the attractive sets. The black line is the mean value over different realizations of the map, while the black dotted lines represent one standard deviation.\\

For $g_1$ observable function, according to equation \ref{g1res}, we expect to find $\xi=0$ in both cases and this is verified by experimental data shown in Figure \ref{c1s1}a). For the scale parameter a similar agreement is achieved in respect to the theoretical parameters $\Delta= \frac{1}{\sigma(g_1)}=\log(2)/\log(3) \simeq  0.6309$ for the Cantor Set and
  $\Delta=\frac{1}{\sigma(g_1)}=\log(3)/ \log(2) \simeq 1.5850$ for Sierpinskij shown in figure \ref{c1s1}b) with a green line.   Eventually, the location parameter $\mu$ shows a logarithm decay with $n$ as expected from equation \ref{g1res}. A linear fit of $\mu$ in respect to $\log(n)$ is shown with a red line in figure \ref{c1s1}c). The linear fit angular coefficient $\kappa$ of equation \ref{g1res} satisfies $\Delta=1/\kappa$ and the dimension computed from data using this relation is $ \Delta = 0.64 \pm 0.01$ for Cantor and  $\Delta = 1.59 \pm 0.01$ in the Sierpinskij triangle.\\
The agreement between the conjecture and the results are confirmed
 even for $g_2$ type and $g_3$ type observable functions shown in figures \ref{c2s2} and \ref{c3s3} respectively. In this case we have experienced some problems in the convergence for the Cantor map has using $\alpha=2$ and $\alpha=3$. The problem is possibly due to the fact that L-moments method works better if $\xi' \in [-0.5, 0.5]$ while for $\alpha\leq3$ the shape parameter $|\xi'|>0.5$ for the Cantor map. For this reason results are shown using $\alpha=4$ both for Sierpinskij triangle and Cantor IFS, for all the experiments presented  the constant value in $g_3$ will be $C=10$. For both observables  $g_2,g_3$ there is strong agreement
between the experimental and  theoretical $\xi'$ values.
In figures \ref{c2s2}b), \ref{c2s2}c), \ref{c3s3}b) a log-log
scale is used to highlight the behavior described by eq. \ref{g2res}
and \ref{g3res}. We can check again the value of the dimension using the angular coefficient $\kappa$ which satisfies $\Delta=1/(\alpha|\kappa|)$. In Table 1 we compare these results with the theoretical values. The error is here represented as one standard
deviation of the ensemble of realizations and we find a good agreement between theoretical and experimental parameters within two standard deviations. 
Eventually, computing $g_3$ observable function we expect to find  a constant value for $\mu(g_3)=C=10$  while $\sigma(g_3)$ has to grow with a power law in respect to $n$ as expected comparing with equation \ref{g3res}.\\
Other tests have been done computing the statistics using  parameter $\alpha=5,6,7,8$ for $g_2$ and $g_3$ observables. Also in this cases, no deviation from the behavior described in eq. \ref{g1res}-\ref{g3res} has been found.

\begin{table}
\begin{center}
 \begin{tabular}{ | l | l| l|}
    \hline
    \textbf{ $\Delta = 1/(\alpha |\kappa|)$} &  \textbf{Cantor} & \textbf{Sierpinskij} \\ \hline
    Theoretical  &  $ \log(2)/\log(3) \simeq  0.6309 $ &  $\log(3)/ \log(2) \simeq 1.5850$   \\ \hline
    $\mu(g_2)$ &  $0.636\pm0.006$ &  $|\kappa|=1.592\pm0.007$  \\ \hline
    $\sigma(g_2)$ &  $0.634\pm0.007$ &  $|\kappa|=1.56\pm0.02$  \\ \hline
    $\sigma(g_3)$ & $0.64\pm0.01$ &  $|\kappa|= 1.62 \pm0.01$  \\ \hline
\end{tabular}
\label{tCS}
\caption{ Information dimension $\Delta$ computed taking the logarithm of equations \ref{g2res}-\ref{g3res} and computing the angular coefficient $\kappa$ of a linear fit of data; IFS with uniform weights and Sierpinskij triangle, $\alpha=4$, $C=10$.}
\end{center}
\end{table}

\subsubsection{IFS with non-uniform weights}

Let us now consider the case of an IFS with different weights:

\begin{equation}
f_k(x)=a_k + \lambda_k x \qquad x \in [0,1] \qquad k=1,2,...,s
\end{equation}

and each $f_i$ is iterated with (different) probability $w_i$.

In this case it is possible to compute the information dimension as
the ratio between the metric entropy and the Lyapunov exponent of
the associated balanced measure  \citep{barnsley2000fractals}: we
get the following expression:

\begin{equation}
\Delta= \frac{w_1\log w1 + ... + w_s \log w_s}{ w_1 \log\lambda_1 + ... + w_s \log \lambda_s}
\label{D1_IFS}
\end{equation}

In this analysis we have considered the following IFS:

\begin{equation}
\begin{cases}
f_1(x)= x/3 & \mbox{ with weight } w \\
f_2(x)= (x+2)/3  & \mbox{ with weight } 1-w \\
\end{cases}
\label{IFS_M}
\end{equation}

and we have changed the weight $w $ between 0.35 and 0.65 with $0.01$ step. For $w=0.5$ we obtain the same results shown in the previous section, while for different
weights we can check the expression \ref{D1_IFS}.\\

The presence of different weights makes the convergence process
sensible to the choice of the sample point $\zeta$ where our
observable reaches its maximum. For that reason we took several
different values of $\zeta$ in order to obtain a reliable estimations
 of the information that
should be obtained, in this case, as an average property.

\begin{figure}[H] \begin{center}
\advance\leftskip-1.5cm

        \includegraphics[width=0.9\textwidth]{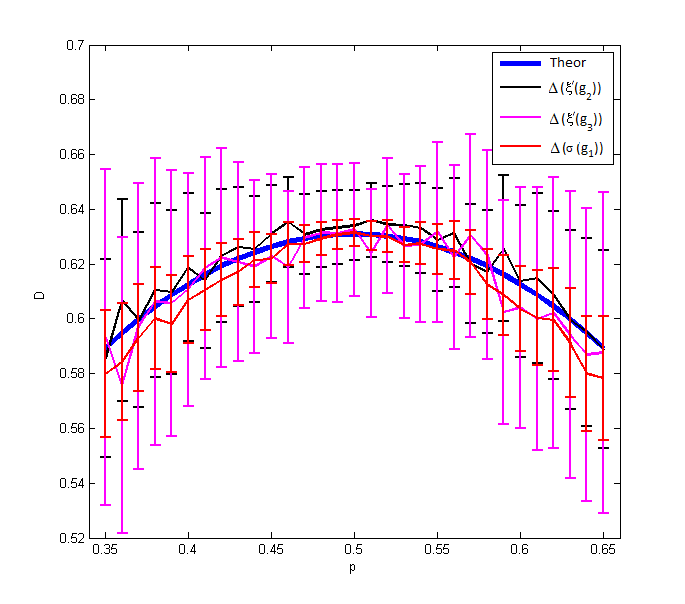}
        \caption{ Different estimation of $\Delta$ dimension obtained using Extreme Value distribution for the IFS in equation \ref{IFS_M}}
      \label{IFS_MF}
    \end{center} \end{figure}

In Figure \ref{IFS_MF} we present the dimension $\Delta$
computed using relationship \ref{g1res}-\ref{g3res}.
 In particular we can compute the dimension from eq. \ref{g1res} as:

\begin{equation}
\Delta(\sigma(g_1))=\frac{1}{<\sigma(g_1)>}
\label{dsigma1}
\end{equation}

We can infer dimension also from  eq.  \ref{g2res}, eq. \ref{g3res} as:

\begin{equation}
\Delta(\xi'(g_i))=\frac{1}{\alpha|<\xi'(g_i)>|}, \quad i=2,3
\label{dcsi}
\end{equation}

and in all expressions above the brackets $<.>$ indicate
 an average on different sample points $\zeta$.
 For the rest of the numerical computations
 we set $\alpha=5$. The parameters have been computed
 using 1000 different initial conditions on the support of the attractor,
 and for 30 realizations of each sample point $\zeta$,the block-maxima approach
 is here used with $n=m=1000$. The error bar are computed using the standard error propagation rules.

The agreement between the theoretical dimension and the experimental data is evident for all the weights and for all the observable considered. The uncertainty increases when $w$ is much different from 0.5. This is due to the fact that as soon as we change the weight to be different from 0.5 the parameters spread increase to take in account the local properties of the attractor.
The best agreement and less uncertainty is
 achieved considering the dimension as computed
  from $\sigma(g_1)$ observable.
  This is possibly due to the slower convergence for $g_2$ or $g_3$
  observables to the respective theoretical distributions:
  in $g_1$ we modulate the distances with a logarithm function
   while in $g_2$ and in $g_3$ power laws are used. Nevertheless, all the data show the right trend.

\subsection{Non-trivial singular measures}

In the previous subsection we have analysed the relatively simple cases of Cantor sets generated with IFS. In order
to provide further support to our conjectures, we now present some application of our theory to the output of dynamical
systems possessing a less trivial singular measures. We consider three relevant examples of two dimensional maps.

\paragraph{The Baker map\\}
The Baker map is defined as follows:

\begin{equation}
x_{t+1}=\begin{cases}
\gamma_a x_t \mod 1 & \mbox{ if } y_t< \alpha \\
1/2 + \gamma_b x_t  \mod 1 & \mbox{ if } y_t \geq \alpha \\
\end{cases}
\label{Ba_x}
\nonumber
\end{equation}

\begin{equation}
y_{t+1}=\begin{cases}
\frac{y_t}{\alpha} \mod 1 &  \quad \qquad \mbox{ if } y_t< \alpha \\
\frac{y_t-\alpha}{1-\alpha}  \mod 1 & \quad  \qquad \mbox{ if } y_t \geq \alpha \\
\end{cases}
\label{Ba_y}
\end{equation}

we consider the classical value for the parameter: $\alpha=1/3, \gamma_a=1/5$ and $\gamma_b=1/4$.

Rigorous analytical  results are available for the estimation
 of the information dimension \citep{kaplan1979chaotic}.
For our parameter values, the analytical expected value is $D\simeq1.4357$.

\begin{figure}[H] \begin{center}
\advance\leftskip-1.5cm
        \includegraphics[width=0.9\textwidth]{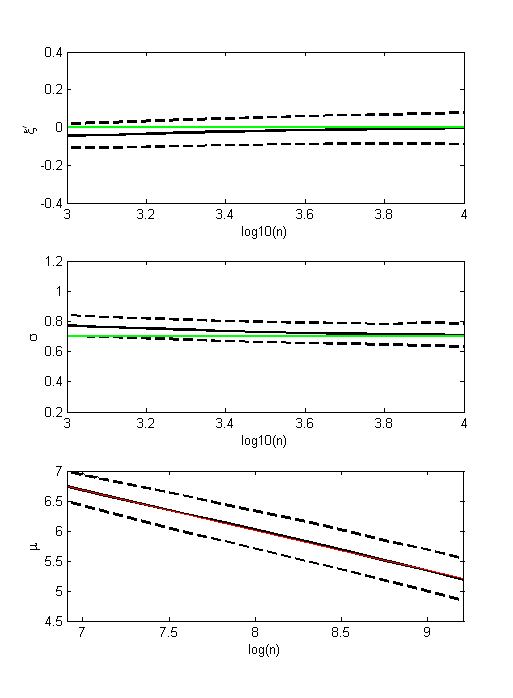}
        \caption{ $g_1$ observable. \textbf{a)} $\xi'$ VS $\log_{10}(n)$; \textbf{b)} $\sigma$ VS $\log_{10}(n)$; \textbf{c)} $\mu$ VS $\log(n)$.  Baker map. Dotted lines represent one standard deviation,  red lines represent a linear fit, green lines are theoretical values.}
      \label{b1}
    \end{center} \end{figure}

\begin{figure}[H] \begin{center}
\advance\leftskip-1.5cm
        \includegraphics[width=0.9\textwidth]{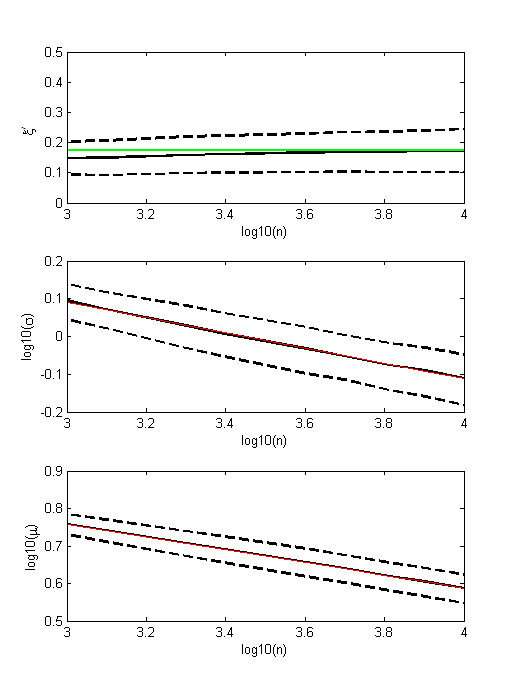}
        \caption{  $g_2$ observable \textbf{a)} $\xi'$ VS $\log_{10}(n)$; \textbf{b)} $\log_{10}(\sigma)$ VS $\log_{10}(n)$; \textbf{c)} $\log_{10}(\mu)$ VS $\log_{10}(n)$. Baker map. Dotted lines represent one standard deviation,,  red lines represent a linear fit, green lines are theoretical values.}
      \label{b2}
    \end{center} \end{figure}

\begin{figure}[H] \begin{center}
\advance\leftskip-1.5cm
        \includegraphics[width=0.9\textwidth]{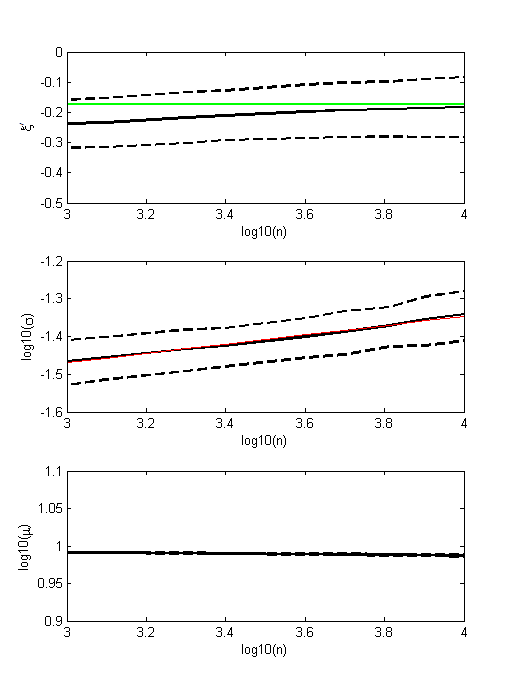}
        \caption{  $g_3$ observable. \textbf{a)} $\xi'$ VS $\log_{10}(n)$; \textbf{b)} $\log_{10}(\sigma)$ VS $\log_{10}(n)$; \textbf{c)} $\log_{10}(\mu)$ VS $\log_{10}(n)$.  Baker map. Dotted lines represent one standard deviation,  red lines represent a linear fit, green lines are theoretical values.}
\label{b3}
\end{center} \end{figure}

We have performed the same analysis detailed in section 4.1, but
with a difference. This map is invertible and its invariant set is
an attractor given by the cartesian product of a segment, along the
$y$-axis, and a one dimensional Cantor set along the $x$-axis. The system possess an invariant SRB measure, which can practically be constructed by taking ergodic sums for any point sitting on the basin of attraction.  In order to compute the center of the balls on the attractor, we proceed in a similar manner as for repellers (see above), namely we take any point  $x$ in the basin of attraction and we iterate it $t$ times with $t$ much bigger than the sequence of observed events.  Then we take $\zeta$ as the point $f^t(x)$: it will be closer and closer to the attractor and distributed according to the SRB measure.  In our set up: $\alpha=4$ for $g_2$ and $g_3$, $C=10$.
The results are shown in figures \ref{b1}-\ref{b3} the black continuous lines will represent the parameter average over different initial conditions and the black dotted lines the standard deviation of the distribution of the estimated parameters.\\
The expected theoretical $\xi'$ values are  within one standard deviation of the results of the fit for all three observable. The agreement seems to be better when we increase $n$ even if this correspond to a decrease of $m$ in our set up. This behavior is quite interesting since it seems that we obtain a much better convergence to theoretical values if $n\simeq 10^4$, while in all other example there is no such a difference between $n=1000$ and $m=10000$. A similar consideration can be made for $\sigma(g_1)$ shown in Figure \ref{b1}b) that approaches the theoretical values value for bigger $n$ values. The angular coefficient of the linear fit for $ \mu(g_1)$ shown in the semilog plot in Figure \ref{b1}c) allow us to estimate the dimension $\Delta=1/|\kappa|=1.48\pm 0.03$ that is consistent with the theoretical values within two standard deviations  .\\

Log-log plots of the parameters against $n$ are shown
in figures \ref{b2}b), \ref{b2}c) and \ref{b3}b), and the value of the dimension $\Delta$  computed using the angular coefficients are reported in Table 2. The agreement with expected value is good enough for all the parameters and better for $\mu(g_2)$.
Eventually, Figure  \ref{b3}c) shows that $\mu(g_3)$ approaches $C=10$.

\newpage

\paragraph{The H\`enon and Lozi maps\\}

The H\`enon map is defined as:

\begin{equation}
\begin{array}{lcl}
x_{t+1}&=&y_t +1 -a x_t^2\\
y_{t+1}&=&b x_t\\
\end{array}
\end{equation}

while in the Lozi map $a x_t^2$ is substituted with $a |x_t|$ :

\begin{equation}
\begin{array}{lcl}
x_{t+1}&=&y_t +1 -a |x_t|\\
y_{t+1}&=&b x_t\\
\end{array}
\end{equation}

We consider the classical set of parameter $a=1.4$, $b=0.3$ for the H\`enon map
and $a=1.7$ and $b=0.5$ for the Lozi map.

\citet{young1985bowen} proved the existence of the SRB measure for the Lozi map, whereas for the 
H\`enon map no such rigorous proof exists, even if convincing numerical results suggest its existence \citep{badii1987renyi}. Note that \citet{benedicks1991dynamics} proved the existence of an SRB measure for the H\`enon map with a different set of parameters. 
Using the classical Young results which makes use of the Lyapunov exponents, we obtain an exact result for $\Delta$ for the Lozi attractor:

$$\Delta \simeq 1.40419 $$

Instead, in the case of the H\`enon attractor, we consider the numerical estimate provided by \citet{grassberger1983generalized}:

$$\Delta=1.25826 \pm 0.00006$$

\begin{figure}[H] \begin{center}
\advance\leftskip-1.5cm
        \includegraphics[width=0.9\textwidth]{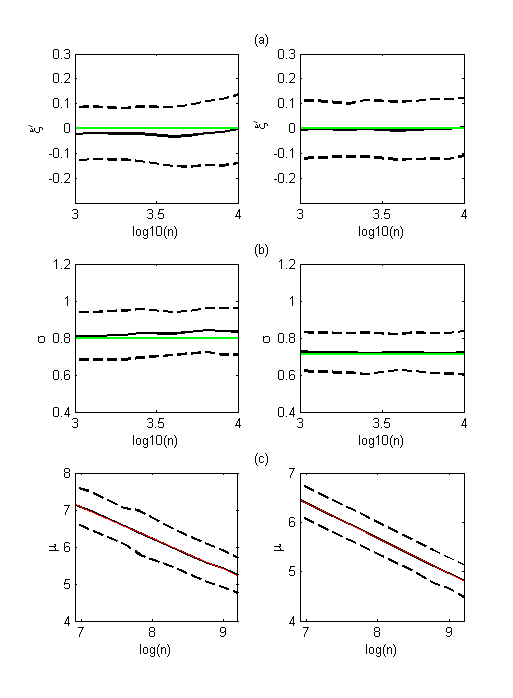}
        \caption{ $g_1$ observable. \textbf{a)} $\xi'$ VS $\log_{10}(n)$; \textbf{b)} $\sigma$ VS $\log_{10}(n)$; \textbf{c)} $\mu$ VS $\log(n)$. Left: H\`enon map, Right: Lozi map. Dotted lines represent one standard deviation,  red lines represent a linear fit, green lines are theoretical values.}
      \label{h1l1}
    \end{center} \end{figure}

\begin{figure}[H] \begin{center}
\advance\leftskip-1.5cm
        \includegraphics[width=0.9\textwidth]{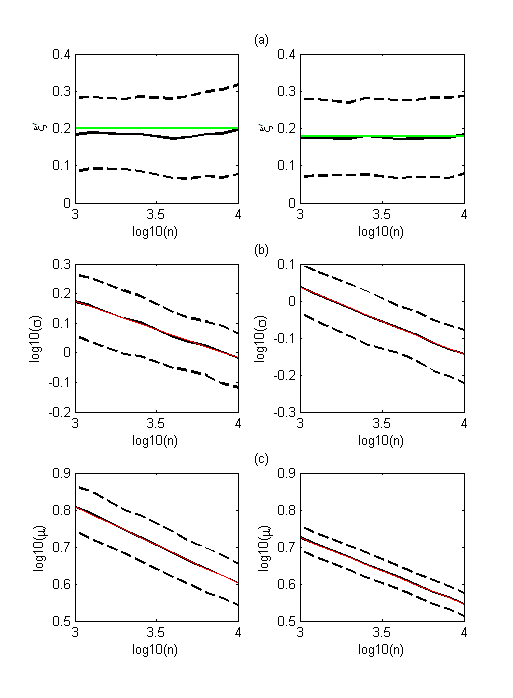}
        \caption{  $g_2$ observable \textbf{a)} $\xi'$ VS $\log_{10}(n)$; \textbf{b)} $\log_{10}(\sigma)$ VS $\log_{10}(n)$; f{c)} $\log_{10}(\mu)$ VS $\log_{10}(n)$. Left: H\`enon map, Right: Lozi map. Dotted lines represent one standard deviation,,  red lines represent a linear fit, green lines are theoretical values.}
      \label{h2l2}
    \end{center} \end{figure}

\begin{figure}[H] \begin{center}
\advance\leftskip-1.5cm
        \includegraphics[width=0.9\textwidth]{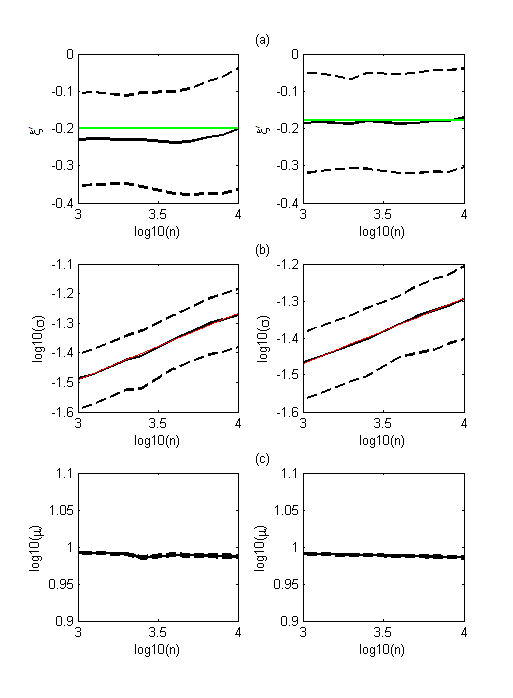}
        \caption{  $g_3$ observable. \textbf{a)} $\xi'$ VS $\log_{10}(n)$; \textbf{b)} $\log_{10}(\sigma)$ VS $\log_{10}(n)$; \textbf{c)} $\log_{10}(\mu)$ VS $\log_{10}(n)$.  Left: H\`enon map, Right: Lozi map. Dotted lines represent one standard deviation,  red lines represent a linear fit, green lines are theoretical values.}
\label{h3l3}
\end{center} \end{figure}

As in the previous cases the GEV distribution is computed with L-moments methods varying $n$ and $m$ and averaging the distribution parameters over 1000 different sample points chosen as described before for the Baker map. Results are presented in figures \ref{h1l1}-\ref{h3l3}, the  plots on the left-hand side refer to the H\`enon map, while on the right-hand side the results refer to the Lozi map.\\
When considering $\xi'$, the numerical results are in agreement  with the theoretical estimates. Nevertheless, the parameters distribution have a rather range spread which indicates a slower convergence towards the expected values in respect to what is observed for the IFS case. The experimental values of $\sigma(g_1)$ , shown in Figure \ref{h1l1}b , approach the theoretical values shown by a green line. The  angular coefficient computed from the semilog plot of $\mu(g_1)$ represented in Figure \ref{h1l1}c) gives us an estimate of the dimension $\Delta=1/|\kappa|= 1.234 \pm 0.015$  for H\`enon and $\Delta=1/|\kappa|= 1.40 \pm 0.01$ for Lozi.\\

The other angular coefficients related to $g_2$ and $g_3$ observables for the plots shown in figures \ref{h2l2}b), \ref{h2l2}c), and \ref{h3l3}b) are presented in Table 2. Within 2 standard deviation they are comparable with the theoretical ones and the best agreement is achieved considering $\mu(g_2)$. The constant $C=10$ is approached quite well  (Figure \ref{h3l3}c)).

\begin{table}[H]
\begin{center}
 \begin{tabular}{ |l| l | l| l|}
    \hline
    \textbf{$\Delta = 1/(\alpha |\kappa|)$} &  \textbf{Baker} & \textbf{H\`enon} & \textbf{Lozi} \\ \hline
    Theor.  &   $ 1.4357$ & $ 1.2582$ &  $1.4042$   \\ \hline
    $\mu(g_2)$ &   $  1.47 \pm 0.02$ & $1.238 \pm0.009$ &  $1.396\pm0.008$  \\ \hline
    $\sigma(g_2)$ &    $1.39 \pm 0.04$ &$1.35\pm0.07$ &  $1.38\pm0.02$  \\ \hline
    $\sigma(g_3)$ &  $ 1.56 \pm0.08$ & $1.15\pm0.07$ &  $ 1.42 \pm0.01$  \\ \hline
\end{tabular}
\label{tBHL}
\caption{Information dimension $\Delta$ computed taking the logarithm of equations \ref{g2res}-\ref{g3res} and computing the angular coefficient $\kappa$ of a linear fit of data; for Baker, H\`enon and Lozi maps.}
\end{center}
\end{table}

The slower convergence for these maps may be related  to the difficulties experienced computing the dimension with all box-counting methods, as shown in \citet{grassberger1983generalized},\citet{badii1987renyi}. In that case it has been proven that the number of points that are required to  cover a fixed fraction of the attractor support diverges faster than the number of boxes itself for this kind of non uniform attractor. In our case the situation is similar since we consider balls around the initial condition $\zeta$. As pointed out, the best result for the dimension is achieved using the parameters provided by $g_1$ observable since the logarithm modulation of the distance exalts proper extrema while weights less possible outliers.

\section{Final Remarks}

Extreme Value Theory is attracting a lot of interest both in terms of extending pure mathematical results and in terms of applications to many fields
of social and natural science. As an example, in geophysical applications
is crucial to have a tool to understand  and forecast climatic extrema and events such as  strong earthquakes and floods.\\
Whereas the classical Extreme Value Theory deals with stochastic processes, many applications demanded to understand whereas it could rigorously
be used to study the outputs of deterministic dynamical systems. 
The mathematical models used to study them present a rich structure and their attracting sets
are very often strange attractors. In such sense is extremely important to develop an extreme value theory
for dynamical system with singular measures. Recently, The existence of extreme value laws for dynamical  systems
preserving an absolutely continuous invariant measure or a singular
continuous invariant measure has been  proven if strong
mixing properties or  exponential hitting time statistics on balls
are  satisfied.\\
In this work we have extended the results presented in
\citet{faranda2011numerical} to the case of dynamical systems with
singular measures. Our main results is that there exist an extreme
value distribution for this kind of systems that is related to the
GEV distribution when observable functions of the distance between
the iterated orbit and the initial conditions are chosen.
The three extreme
value type for the limit distribution laws for maxima and the
generalized distribution family (GEV) are absolutely continuous
function. We will recover the GEV using histograms on the frequency
of maxima; in this way the cumulative distribution function which we
got from such an histogram will have {\em plateaux} just in
correspondence of the holes of the Cantor set, whenever this one is
the invariant set. This could be easily explained by the very nature
of our observables which measure the distance with respect to a
given point: there will be distances which are not allowed when such
distances are computed from points in the holes. It should be
stressed that such a cumulative distribution function, which is a
sort of devil staircase and therefore is a singular continuous
function, in any way could converge to a GEV distribution.  The
latter as the three type extreme values laws are
 normalized laws which must be adjusted in order to give a
non-degenerate asymptotic distribution. The strength of our
approach, as we said above, relies in the possibility to infer the
nature and the value of such normalizing constants by a fitting
procedure on the unnormalized data, the histograms. This worked very
well for probability measure which were absolutely continuous. We
will see that it works also for singular measures (and the
normalizing constants will be related to the information dimension),
provided we remind that this time the fitting procedure will contain
a sort of extrapolation to smooth out the gaps of the Cantor sets.\\
It is interesting to observe that on
  Cantor sets the notion of {\em generic point} is not so obvious as
   for smooth manifolds which support Lebesgue measure: in this  case
   in fact one could suppose that each point accessible for  numerical
   iterations is generic with respect to an invariant measure which is
    in turn absolutely continuous. This notion of genericity is restored
    on attractor by considering the SRB measure. Instead for Iterated Function Systems we can dispose of uncountably many measures, but we have a precise manner to  identify them and this will be reflected in the different dimensions produced by the numerical computation of the parameters of the GEV. The possibility to discriminate among different singular measures having all the same topological support is another indication of the validity and of the efficiency  of our approach.\\
We have also shown that the parameters of the distribution are intimately related to the
information dimension of the invariant set. We have tested our
conjecture with numerical experiments on different low dimensional
maps such as : the middle third Cantor set, the Sierpinskij
triangle), Iterated Function System (IFS) with non-uniform weights,
strange attractors such as Lozi and H\`enon. In all cases considered
there is agreement between the theoretical parameters and the
experimental ones. The extimates of $\Delta$ are in agreement with the theoretical values in all cases
considered. It is interesting to observe that  the algorithm
described with the selection of maxima acts like a magnifying glass
on the neighborhood of the initial condition.
 In this way we have both a powerful tool to study and highlight the fine structure of the attractor, but, on the other
  hand we can  obtain global properties averaging on different initial conditions.
Even if we are dealing with very simple maps for which many properties are known it is clear from numerical
experiments that is not so obvious to observe a good convergence to the GEV distribution. Even if we are able to compute very
large statistics and the results are consistent with theoretical values, the error range is wide if compared to the experiment
for maps with a.c.i.m. measures that we have carried out in \citet{faranda2011numerical}. This should be taken in consideration
 each time this statistics is applied in a predictive way to more complicated systems.\\
In the case of an experimental temporal series, for which the underlying dynamics is unknown, a classical problem is to obtain the
dimensionality of the attractor of the dynamical systems which generated it. This can be achieved through the so called Ruelle-Takens
 delay embedding where, starting  from the time series of an observable $O(n)$, we can construct the multivariate vectors in a $\Delta$-dimensional
 space:

$$ \phi(N)=[O(n), O(n+1),..., O(n+d-1)] $$

and study the geometrical properties using the Recurrence Qualification Analysis \citep{marwan2007recurrence}. The minimum value of $\Delta$ needed to reconstruct the actual dimension $d^*$ is given by  $[2d]+2$. Using the procedure described in this paper we could find an estimate for $d^*$ thus determining the minimum value of $\Delta$ to be used in the Recurrence Quantification Analysis. We will test this strategy in a subsequent publication.

The theory and the algorithm presented in this work and in \citet{faranda2011numerical} allow to study in detail the recurrence of an orbit around a point:  this is due to the particular choice of the observables that require to compute distances between initial and future states of the system. Understanding the behavior of a dynamical system in a neighborhood of a particular initial condition is of great interest in many applications. As an example, in weather forecast and climate it is important to  study the recurrence of patterns (the so called analogues). In principle, applying the  extreme value statistics to the output of meteorological models, will make possible to infer dynamical properties related to the closest return towards a certain weather pattern. EVT will give information not only about the probability distribution of the extrema but also about the  scaling of the measure of a ball centered on the chosen initial condition  providing an insight to the dynamical structure of the attractor.

\section{Acknowledgments}

S.V. was supported by the CNRS-PEPS Project {\em Mathematical Methods of Climate Models }, and he thanks the GDRE Grefi-Mefi for having supported exchanges with Italy.
V.L. and D.F. acknowledge the financial support of the EU FP7-ERC project NAMASTE ``Thermodynamics of the Climate System". The authors are grateful to Jorge Freitas for many helpful comments.

\bibliographystyle{plainnat}
\bibliography{Cantor_GEV}

\end{document}